\input amstex
\documentstyle{amsppt}
\NoRunningHeads
\magnification=\magstep1


\def\a{\alpha}
\def\b{\beta}
\def\g{\gamma}
\def\c{\gamma}
\def\s{{\sigma}}
\def\cO{\Cal O}

\def\th{\theta}
\def\t{\theta}

\def\tS{{\tilde S}}
\def\-{{\setminus}}
\def\<{{\langle}}
\def\>{{\rangle}}

\def\defi{{def}}
\def\lra{{\longrightarrow}}
\def\GL{{\hbox{GL}}}

\def\NO{\negmedspace\negmedspace\negmedspace\negmedspace\negmedspace
\negmedspace\negmedspace\negmedspace\negmedspace\negmedspace\negmedspace
\negmedspace\negmedspace\negmedspace\negmedspace\negmedspace\negmedspace
\negmedspace\negmedspace\negmedspace\negmedspace\negmedspace\negmedspace
\negmedspace\negmedspace\negmedspace\negmedspace\negmedspace\negmedspace
------}

\topmatter
\title
Commutation relations among 
\\
quantum minors 
in $\cO_q(M_n(k))$
\endtitle

\author R. Fioresi
\endauthor
\thanks Investigation supported by
the University of Bologna, funds for selected research topics. \endthanks
\affil
Dipartimento di Matematica, Universita' di Bologna
\endaffil
\address
Piazza Porta San Donato 5, 40126 Bologna, Italy
\endaddress
\email
rfioresi\@math.ucla.edu
\endemail

\abstract  In this paper we compute explicit
formulas for the commutation relations between any two
quantum minors in the quantum matrix bialgebra
$\cO_q(M_n(k))$. The product of any two
minors is expressed as linear combination
of products of minors strictly less in certain 
orderings. 
\endabstract

\endtopmatter

\document
\bigskip

{\bf 1. Introduction}

\medskip

The problem of computing commutation relations among the quantum
minors in the quantum matrix bialgebra has been extensively studied
especially in connection with the construction of geometric 
homogeneous quantum spaces. Some of the first computations
of commutation relations appear 
in [LR] using the $R$-matrix formalism and later in [TT]
using the matrix bialgebra notation.
In both cases the purpose is to give the
quantum deformation of the coordinate ring of the quantum flag
manifold. Such ring is in fact generated inside
the quantum matrix ring by certain quantum minors.
However in both papers no explicit 
formula is given for the commutation
among the quantum minors considered.
Explicit calculations for some of the commutation rules
among quantum minors appear also in [Do1], [Do2] where the
purpose instead is to study the $q$-difference intertwining
operator.
\medskip

An explicit formula for the commutation relation of all commutation
relations of two minors
having the same size and same columns appears in [Fi1]. 
Later in [Fi2] is considered also the case of commutation relations
among quantum minors with possibly different size, but with the
first columns in common. In both cases the purpose of the formula
is to obtain a quantum deformation of the grassmannian
and flag manifolds.
Similar formulas appear also in [GL], with a different purpose.
They are used to prove that the matrix bialgebra ring admits a 
basis made by certain products of quantum minors, as it happens for the
commutative case.
\medskip

The present paper has the purpose to compute the most general
commutation relations among any two quantum minors chosen
in the quantum matrix bialgebra.
In the relation the product of two minors will be expressed as
linear combination of products of minors strictly less in certain
specified orderings.
\medskip

This paper is organized as follows.

\medskip

In \S 2 we state and discuss the Manin relations, that is the
relations that define the quantum matrix bialgebra in a free
tensor algebra. In particular we describe two sets of generators for
the ideal of the Manin relations $I_M$:
the $R$ relations and the $S$ relations. 
The $R$ relations (respectively the $S$ relations) allow us to
rewrite a product of two elements in terms of a sum of products in
which the row  order (respectively the column order) is preserved.
\medskip
In \S 3 we define the type $F$ and type $G$ operations. These operations
amount to perform a rewriting 
of an certain expression in a free algebra
using the Manin relations. The operations of type $F$ make use
of the $R$ relations, the operations of type $G$ of the $S$ relations.
These operations appeared in less general
context in [Fi1]. 
\medskip

In \S 4 we define quantum row and column minors
with rows $I$ and columns $J$: these are elements
in the free algebra $k_q\<a_{ij}\>$. They both
specialize to the quantum minor with rows $I$ and columns $J$ in the
matrix bialgebra $\cO_q(M_n(k))=k_q\<a_{ij}\>/I_M$.
The operations of type $F$ will act on products of
row minors, the operations of type $G$ on products of column minors.
We also define the notion of commutation relation and of
intertwining row and column order of an index. 

\medskip

In the last sections we derive the formulas for the commutation
relation of any two minors in $\cO_q(M_n)k))$. For clarity reasons we
prefer to write the formulas in special cases,
that cover all the possibilities. A list of all the cases and
the theorems appears at the end of \S 7.

\medskip

Acknoledgements. We wish to thank Prof. T. 
Lenagan for many helpful discussions
and comments. 

\bigskip

{\bf 2. Preliminaries: The Manin relations.}
\medskip
Let $k$ be a field.

Let $\cO_q(M_n(k))$ denote the quantum matrix bialgebra.
This is the algebra generated over $k_q=k[q,q^{-1}]$ ($q$ being
an indeterminate) by the non commutative 
indeterminates $x_{ij}$, $i,j=1 \dots n$, subject to the
relations:
$$
\matrix
x_{ij}x_{ik}=q^{-1}x_{ik}x_{ij}, \qquad
x_{ji}x_{ki}=q^{-1}x_{ki}x_{ji},\quad j<k \cr\cr
x_{ij}x_{kl}=x_{kl}x_{ij}, \quad i<k,j>l \hbox{ or } i>k,j<l \cr\cr
x_{ij}x_{kl}-x_{kl}x_{ij}=(q^{-1}-q)x_{kj}x_{il}, \quad i<k,j<l
\endmatrix
$$
These are called the {\it Manin relations}. 
\medskip

We want to rewrite the Manin relations in a way more suitable  
to our purpose, that is the commutation relations between quantum minors.
\medskip
Let $k_q\<a_{rs}\>$ denote the free tensor algebra over $k_q$ generated
by the indeterminates $a_{rs}$, $r,s = 1 \dots n$. 

Let $I_M$ denote the ideal of
the Manin relations in $k_q\<a_{rs}\>$. All congruences have now to
be intended taking place in $k_q\<a_{rs}\>$ modulo $I_M$ unless otherwise
specified.

\medskip

{\bf Definition (2.1)}.
Let's define the {\it relation $R$}, 
$R_{a_{ij}a_{kl}}(a_{ij}a_{kl}) \in k_q\<a_{rs}\>$,
commuting the elements $a_{ij}$, $a_{kl}$ as:

$$
R_{a_{ij}a_{kl}}(a_{ij}a_{kl})
=a_{ij}a_{kl}-q^{A(i,j,k,l)}q^{B(i,j,k,l)}a_{kl}a_{ij}
-C(i,j,k,l)(q^{-1}-q)a_{kj}a_{il} 
$$

where the coefficients $A$, $B$, $C$ are:

$$
\matrix
A(i,j,k,l)=\cases -1 & $if $ i<k, j=l \cr 
                  +1 & $if $ i>k, j=l \cr 
                   0 & $ otherwise $\endcases
\cr \cr \cr
B(i,j,k,l)=\cases  
                  -1 & $if $j<l, i=k \cr 
                  +1 & $if $j>l, i=k \cr 
                   0 & $ otherwise $\endcases
\cr \cr \cr
C(i,j,k,l)=\cases
                  +1 & $ if $i<k, j<l \cr
                  -1 & $ if $i>k, j>l \cr
                  0 & $ otherwise $ \endcases
\endmatrix
$$

So for example the commutation relation of $a_{11}a_{22}$ is

$$
R_{a_{11}a_{22}}=a_{11}a_{22}-a_{22}a_{11}-(q^{-1}-q)a_{21}a_{12}
$$

since $A(1122)=0$, $B(1122)=0$, $C(1122)=1$.

Notice that $R_{a_{ij}a_{kl}}(a_{ij}a_{kl}) \equiv 0$ modulo $I_M$.

\medskip

We also define
the {\it relation $S$},
$S_{a_{ij}a_{kl}}(a_{ij}a_{kl}) \in k_q\<a_{rs}\>$ commuting elements 
$a_{ij}$, $a_{kl}$ as:
$$
\matrix
S_{a_{ij}a_{kl}}(a_{ij}a_{kl})
=a_{ij}a_{kl}-q^{A(i,j,k,l)}q^{B(i,j,k,l)}a_{kl}a_{ij}
-C(i,j,k,l)(q^{-1}-q)a_{il}a_{kj} 
\endmatrix
$$
where the coefficients $A$, $B$, $C$ have been defined above.

For example the $S$ commutation relation of $a_{11}a_{22}$ is
$$
S_{a_{11}a_{22}}=a_{11}a_{22}-a_{22}a_{11}-(q^{-1}-q)a_{12}a_{21}
$$
\medskip
{\bf Remark (2.2)}.
Notice that the only difference between the relations $R$'s 
and $S$'s, is in the $(q^{-1}-q)$ term. For $R_{a_{ij}a_{kl}}(a_{ij}a_{kl})$ 
it is $C(i,j,k,l)(q^{-1}-q)a_{kj}a_{il}$
for $S_{a_{ij}a_{kl}}(a_{ij}a_{kl})$ it is $C(i,j,k,l)(q^{-1}-q)a_{il}a_{kj}$. 
It is important to notice that while in $R_{a_{ij}a_{kl}}(a_{ij}a_{kl})$
the last two terms have the {\it row} indices in the same order i.e. $(k,i)$
in $S_{a_{ij}a_{kl}}(a_{ij}a_{kl})$  
the last two terms have {\it column} indices in the same order i.e. $(l,j)$.
This observation is crucial when we will make use of such relations.

The two sets of relations $R_{a_{ij}a_{kl}}(a_{ij}a_{kl})$ 
and $S_{a_{ij}a_{kl}}(a_{ij}a_{kl})$ generate in $k_q\<a_{rs}\>$
the same two sided ideal, namely $I_M$ the ideal of the Manin relations.

\medskip

{\bf Definition (2.3)}. We define {\it commutation relation $R$ of
$a_{ij}a_{kl}$ applied to $a_{\a\beta}a_{\g\delta}$} as:
\medskip

$$
\align
R_{a_{ij}a_{kl}}(a_{\a\beta}a_{\g\delta}) & = 
a_{\a \beta}a_{\c \delta}-q^{A(i,j,k,l)}q^{B(i,j,k,l)}a_{\c \delta}a_{\a \beta}
\cr\cr
 & -C(i,j,k,l)(q^{-1}-q)a_{\c \beta}a_{\a \delta} 
\endalign
$$

Similarly we define {\it commutation relation $S$ of
$a_{ij}a_{kl}$ applied to $a_{\a\beta}a_{\g\delta}$} as:
$$
\align
S_{a_{ij}a_{kl}}(a_{\a\beta}a_{\g\delta}) & =
a_{\a \beta}a_{\c \delta}-q^{A(i,j,k,l)}
q^{B(i,j,k,l)}a_{\c \delta}a_{\a \beta}+
\cr\cr
& -C(i,j,k,l)(q^{-1}-q)a_{\a \delta}a_{\c \beta}
\endalign
$$

Notice that in general:
$R_{a_{ij}a_{kl}}(a_{\alpha \beta} a_{\gamma\delta})$, 
$S_{a_{ij}a_{kl}}(a_{\alpha \beta} a_{\gamma\delta}) \not\in I_M$.
\medskip

The key to the commutation relations between
quantum minors is the following lemma. This lemma
tells us what is the ``error'' when we commute two elements
$a_{ij}a_{kl}$ using the commutation relation $R$ (respectively $S$) 
associated to $a_{kj}a_{il}$ (respectively $a_{il}a_{kj}$).

\medskip

{\bf Lemma (2.4)}.
$$
\matrix
1. \quad R_{a_{ij}a_{kl}}(a_{ij}a_{kl}) = \cases
R_{a_{kj}a_{il}}(a_{ij}a_{kl}) +
(q^{-1} -q)a_{kj}a_{il} & $ if $i>k \cr\cr
R_{a_{kj}a_{il}}(a_{ij}a_{kl}) -
(q^{-1} -q)a_{kj}a_{il} & $ if  $i<k
\endcases
\cr\cr\cr
2. \quad S_{a_{ij}a_{kl}}(a_{ij}a_{kl}) = \cases
S_{a_{il}a_{kj}}(a_{ij}a_{kl}) +
(q^{-1} -q)a_{il}a_{kj} & $ if  $j>l \cr\cr
S_{a_{il}a_{kj}}(a_{ij}a_{kl}) -
(q^{-1} -q)a_{il}a_{kj} & $ if  $j<l
\endcases
\endmatrix
$$

\medskip
{\bf Proof}. (1) is in [Fi1] pg 423. (2) is the same.

\medskip

Let's see an example.
$$
R_{a_{11}a_{22}}(a_{11}a_{22})=R_{a_{21}a_{12}}(a_{11}a_{22})-
(q^{-1}-q)a_{21}a_{12}
$$
In fact $R_{a_{21}a_{12}}(a_{11}a_{22})=a_{11}a_{22}-a_{22}a_{11}$.


\bigskip
{\bf 3. The operations of type $F$ and $G$}.

\medskip
 
The operations of type $F$ and of type $G$
transform a tensor $s \in k_q\<a_{ij}\>$
into another tensor $s'$ such that 
$s \equiv s'$ modulo $I_M$, the ideal of
the Manin relations. These operations allow us to take control on
how minors commute.

\medskip
Let $k_q\<a_{rs}\>_p$ denote the homogeneous component of degree $p$ in 
$k_q\<a_{rs}\>$. Let $n_0$ be a fixed integer.

\medskip
Let $T \subset k_q\<a_{rs}\>_{n_0}$ 
be the free $k_q$-module in $k_q\<a_{rs}\>$ generated by the 
monomial tensors:
$$
\t_{PQ}=a_{p_1q_1} \dots a_{p_{n_0}q_{n_0}}
$$
with $(p_1 \dots p_{n_0})$, $(q_1 \dots q_{n_0}) \in S_{n_0}$ the permutation group, where we interpret an index $(x_1 \dots x_y)$ with distinct elements
as the permutation $\s \in S_y$, $\s(i)=x_i$.

\medskip

Let's see an example on how these monomial tensors arise and why they
are important for the commutation of quantum minors.

\medskip
{\bf Example (3.1)}.
Consider the following quantum minors (no
row and no column indices in common) in $\cO_q(M_n(k))$:
$$
\matrix
[12,24]=_{\defi}x_{12}x_{24}-q^{-1}x_{22}x_{14}, &
[34,13]=_{\defi}x_{31}a_{43}-q^{-1}x_{33}x_{41}
\endmatrix
$$
Let $p$ be the projection $k_q\<a_{ij}\> \lra 
\cO_q(M_n(k))\cong k_q\<a_{ij}\>/I_M$. We are interested in two
distinct elements 
$s, t \in k_q\<a_{ij}\>$ that project into the product
$[34,13][12,24] \in \cO_q(M_n(k))$:
$$
\align
s & = a_{31}a_{43}a_{12}a_{24}
-q^{-1}a_{31}a_{43}a_{14}a_{22}-q^{-1}a_{33}a_{41}a_{12}a_{24}
+q^{-2}a_{33}a_{41}a_{14}a_{22}
\cr\cr
t & = a_{31}a_{43}a_{12}a_{24}
-q^{-1}a_{31}a_{43}a_{22}a_{14}-q^{-1}a_{41}a_{33}a_{12}a_{24}
+q^{-2}a_{41}a_{33}a_{22}a_{14}
\endalign
$$
Clearly $s \equiv t$. 
Both $s$ and $t$ are sums of monomial tensors, hence $s$, $t \in T$. 
It is important to notice than in $s=\sum \t_{PQ}$, 
the monomial tensors $\t_{PQ}$
appear all with the same row index $P= (3412)$. 
Similarly in $t$
they appear all with the same column index $Q=(1324)$.
This property is shared by tensors corresponding to products
higher order minors and depends on whether
we are choosing row or column expansion to write the minors.

\medskip

Commuting the two minors $[34,13]$, $[12,24]$ corresponds in the
free algebra $k_q\<a_{ij}\>$ to bring all monomials in $s$ to have
row index $(1234)$, using elements in $I_M$. 
This will be achieved using the commutation relations
$R$. Equivalently one could also bring all monomials in $t$ to
have column order $(2413)$. This will be obtained with the commutation
relations $S$. Both these procedures that we are about to describe 
in a more general setting,
will produce the same result in  $\cO_q(M_n(k))$.
In [Fi1] \S 2, is described how to proceed in the first case,
assuming there is no permutation on the column indices. We need
a generalization of the arguments appearing there.

\medskip

We now give the definition of operations of type $F$ and $G$
applied to a monomial tensor. This is a central notion for the commutation
of the minors. Applying an operation of type $F$ (or of type $G$) to 
the monomial tensor $\t$ basically consists in replacing an occurrence of
$a_{ij}a_{kl}$ in $\t$ with a suitable expression, 
so that the result is congruent to $\t$ modulo $I_M$.

\medskip

{\bf Definition (3.2)}. 
Given $\t \in {T}$ a monomial tensor, 
define $F_{a_{ij}a_{kl}} \t $ {\it an operation of type $F$ on $\t$}
as the element of $T$ obtained from 
$\t$ by replacing $a_{ij}a_{kl}$  with:
$$
a_{kl}a_{ij}
+C(i,j,k,l)(q^{-1}-q)a_{kj}a_{il}
$$
If we have no occurrence of $a_{ij}a_{kl}$ in $\t$
we set $F_{a_{ij}a_{kl}} \t =\t$.

We also define $G_{a_{ij}a_{kl}} \t $ {\it an operation of type $G$ on $\t$}
as the element of $T$ obtained from 
$\t$ by replacing $a_{ij}a_{kl}$  with:
$$
a_{kl}a_{ij}
+C(i,j,k,l)(q^{-1}-q)a_{il}a_{kj}
$$
As before if $a_{ij}a_{kl}$ does not occur
in $\t$ we set $G_{a_{ij}a_{kl}} \t =\t$.

If $t= \sum b_{PQ} {\t}_{PQ} \in T$, $F_{a_{ij}a_{kl}}t$ (respectively
$G_{a_{ij}a_{kl}}t$)
is obtained by performing the operation $F_{a_{ij}a_{kl}}$
(respectively $G_{a_{ij}a_{kl}}$) on only one
of the monomial tensors in $t$ and of course we need to specify which one.

We call {\it composition} of two operations of type $F$ (respectively $G$), 
$F_1$ and $F_2$
applied to $t \in T$ the expression: $F_1(F_2t)$.
If two different operations of type $F$ (respectively $G$)
performed on $t \in T$,  
act on different monomial
tensors in $t$, they obviously commute, i.e. it does not matter
which one is performed first on $t$.

\medskip
Notice that there is no ambiguity in the definition since in each
monomial tensor there will be at most one occurrence of $a_{ij}a_{kl}$.
\medskip

{\bf Proposition (3.3)}. {\it Let $\t \in T$. Then:

\noindent 
1. $F_{a_{ij}a_{kl}}\t \equiv \t$ modulo $I_M$.

\noindent
2. $G_{a_{ij}a_{kl}}\t \equiv \t$ modulo $I_M$.
}

{\bf Proof.} Immediate from the definition of operation of type
$F$ and $G$.
\medskip

{\bf Example (3.4)}.
Consider the monomial tensor
$\t_{(3412),(1324)}=a_{31}a_{43}a_{12}a_{24}$:
$$
G_{a_{43}a_{12}}a_{31}a_{43}a_{12}a_{24}=
a_{31}a_{12}a_{43}a_{24}-(q^{-1}-q)a_{31}a_{42}a_{13}a_{24}
$$
We have passed from a column index $Q=(1324)$ 
in $\t_{(3412),(1324)}$, 
to a column index $\tilde Q=(1234)$ in both monomial
tensors of $G_{a_{43}a_{12}}a_{31}a_{43}a_{12}a_{24}$. 

\medskip

In general the result of the operations of type $G$ on a generic $\t_{PQ}$ 
consists of monomial tensors $\t_{P_i \tilde Q}$,
having all the {\it same} column index $\tilde Q$ that is 
a permutation of $Q$.

The same is true for operation of type $F$. The result of such operation
on a monomial tensor $\t_{PQ}$ consists of
monomial tensors $\t_{\tilde PQ_j}$,
having all the {\it same} row index $\tilde P$ that is
a permutation of $P$.
\medskip

We now want to define two natural actions of the group $S_{n_0}$ on 
$T=k_q\<a_{ij}\>_{n_0}$.

\medskip
Given a multiindex $(x_1 \dots x_p)$, $1 \leq x_1 \dots x_p \leq n$ 
define 
$$
\s(x_1 \dots x_p)=_{def}(\s(x_1) \dots \s(x_p)), \qquad \s \in S_p.
$$
\medskip
{\bf Definition (3.5)} Let $t=\sum \t_{P_i,Q_i} \in T$.
Define a {\it row action} of $S_{n_0}$ on $T$ as:
$$
\s t=\sum \t_{\s(P_i),Q_i}, \qquad \s \in S_{n_0}.
$$
Define a {\it column action} of $S_{n_0}$ on $T$ as:
$$
\s t=\sum \t_{P_i,\s(Q_i)}, \qquad \s \in S_{n_0}.
$$

\medskip

We now state a fundamental lemma.

\medskip
{\bf Lemma (3.6)}. {\it Let $\t \in {T}$ be a monomial tensor,
$\sigma=(m,m+1) \in S_{n_0}$. 

Then
\hfil\break
i) If $\s \in S_{n_0}$ acts with a row action:
$$
F_{a_{\sigma(i)j}a_{\sigma(k)l}}\sigma(\t) = \cases
\sigma(F_{a_{ij}a_{kl}}\t ) & $ if $ \{i,k\} \neq \{m,m+1\} \cr\cr
\sigma(F_{a_{ij}a_{kl}}\t ) + (q^{-1}-q)\t &$ if $(i,k) = (m+1,m) \cr\cr
\sigma(F_{a_{ij}a_{kl}}\t ) - (q^{-1}-q)\t &$if $(i,k) = (m,m+1)
\endcases
$$ 

ii) If $\s \in S_{n_0}$ acts with a column action:
$$
G_{a_{i\sigma(j)}a_{k\sigma(l)}}\sigma(\t) = \cases
\sigma(G_{a_{ij}a_{kl}}\t ) & $ if $ \{j,l\} \neq \{m,m+1\} \cr\cr
\sigma(G_{a_{ij}a_{kl}}\t ) + (q^{-1}-q)\t &$ if $(j,l) = (m+1,m) \cr\cr
\sigma(G_{a_{ij}a_{kl}}\t ) - (q^{-1}-q)\t &$if $(j,l) = (m,m+1)
\endcases
$$ 
}

{\bf Proof.} This is a consequence of Lemma (2.4).
\medskip

{\bf Example (3.7)}.
$\t=a_{31}a_{42}a_{13}a_{24}$, $\s=(2,3)$, $j=2$, $l=3$, $\s$
acts with a column action.
$$
\matrix
G_{a_{43}a_{12}}a_{31}a_{43}a_{12}a_{24}=
a_{31}a_{12}a_{43}a_{24}-(q^{-1}-q)a_{31}a_{42}a_{13}a_{24}. \cr\cr
\sigma(G_{a_{42}a_{13}}a_{31}a_{42}a_{13}a_{24})
-(q^{-1}-q)a_{31}a_{42}a_{13}a_{24}= 
\cr\cr
\s(a_{31}a_{13}a_{42}a_{24})-
(q^{-1}-q)a_{31}a_{42}a_{13}a_{24}
\endmatrix
$$

\medskip

{\bf Remark (3.8).} Observe that the operations of type $F$ (respectively G)
can be applied also to monomial tensors with repeated column
(respectively row) indices. The Lemma (3.6) will still be
true with no modifications, in fact the proof of part (i) (respectively (ii))
involves only reasoning on row (respectively column) indices.

\medskip
 
We now define an ordering on the row and column indices
of tensor monomials. 

\medskip

{\bf Definition (3.9)}. Let $K=(k_1 \dots k_r)$,
$1 \leq k_1 <\dots <k_r \leq n_0$
and let $K'$ be its complement in $(1 \dots n_0)$. Both indices
are assumed to be ordered.

We say that a monomial tensor $\t_{PQ}$
is in {\it row order $(K,K')$} if
$P=(K,K')$. A tensor $t =\sum b_i \t_{P_iQ_i} \in T$
is in row order $(K,K')$ if each monomial
tensor $\t_{P_iQ_i}$ 
is in row order $(K,K')$, i.e. $P_i=(K,K')$ for all $i$.
 
Similarly for columns. If
$L=(l_1 \dots l_r)$ 
$1 \leq l_1 <\dots l_r \leq n$
and $L'$ is its complement in $(1 \dots n)$, we say that $\t_{PQ}$ 
is in {\it column order $(L,L')$} if 
$Q=(L,L')$. A tensor $t =\sum b_i \t_{P_iQ_i} \in T$ 
is in column order $(L,L')$ if $Q_i=(L,L')$ for all $i$. 

\medskip
\medskip

We now want to make sure 
that given a generic monomial tensor it is always possible with
a composition of operations of type $F$ (respectively of type $G$)
to obtain a tensor in a given row order (respectively column order).
Clearly the resulting 
tensor will in general be different from the starting monomial
tensor, however they will be congruent modulo $I_M$.

\medskip

{\bf Lemma (3.10)}. {\it Let ${\t}_{PQ}\in T$, $K$ a row index,  
$L$ an column index as in Definition (3.9). 

i)  There exists  $F^K_{{\t}_{PQ}}$, a composition of operations of 
type $F$, such that $F^K_{{\t}_{PQ}}{\t}_{PQ}$ is in order $(K,K')$.
Moreover $F^K_{{\t}_{PQ}}$ can be chosen such that
given two indices $p_i,p_j \in P$, $i<j$, 
the two corresponding elements in each monomial tensor 
of $F^K_{{\t}_{PQ}}{\t}_{PQ}$
have been interchanged once by some operation of type
$F$ in $F^K_{{\t}_{PQ}}$ if 
$K=(\dots p_j \dots p_i \dots)$
and have not been interchanged 
if $K=(\dots p_i \dots p_j \dots)$.

ii) There exists  $G^L_{{\t}_{PQ}}$, a composition of operations of 
type $G$, such that $G^L_{{\t}_{PQ}}{\t}_{PQ}$ is in order $(L,L')$.
Moreover $G^L_{{\t}_{PQ}}$ can be chosen such that
given two indices $q_i,q_j \in Q$, $i<j$, 
the two corresponding elements in each monomial tensor 
of $G^L_{{\t}_{PQ}}{\t}_{PQ}$
have been interchanged once by some operation of type
$G$ in $G^L_{{\t}_{PQ}}$ if $L=(\dots q_j \dots q_i \dots)$
and have not been interchanged 
if $L=(\dots q_i \dots q_j \dots)$. 
}

\medskip

{\bf Proof}. (i) is proven in [Fi1] in Lemma (2.17). (ii) is
the same.

\medskip
We denote with $F^K_t$ (respectively $G^L_t$) the sequence of
operations of type $F$ (respectively $G$) such that
$F^K_tt$ (respectively $G^L_tt$) is in row order $(K,K')$ 
(column order $(L,L')$). The existence of such $F^K_t$ and
$G^L_t$ is granted by Lemma (3.10).
\medskip

We are ready for the lemma which is the corner stone on  which the formula
for the commutation relations is based. 

\medskip

{\bf Lemma (3.11)}. {\it 
Let $F^K_{PQ}$, $G^L_{PQ}$ and  $\th_{PQ}$ as above.
Let $\s=(m,m+1) \in S_{n_0}$.

1. Row version. Assume $\s$ acts on rows.
Then
$$
F^{\s(K)}_{\theta_{\s(P)Q}}\theta_{\s(P)Q}=
\s(F^{K}_{\theta_{PQ}}\theta_{PQ})+(q^{-1}-q)E(m,K,P)\th_0
$$
where $\th_0 \equiv \theta_{PQ}$ and
$$
E(m,K,P)=\cases 0 & $ if $(K,K')=(\dots m \dots m+1 \dots),
P=(\dots m \dots m+1 \dots) \cr
& $or $(K,K')=(\dots m+1  \dots m \dots),
P=(\dots m+1 \dots m \dots) \cr\cr
1 & $ if  $(K,K')=(\dots m \dots m+1 \dots),
P=(\dots m+1 \dots m \dots) \cr\cr
-1 & $ if $(K,K')=(\dots m+1 \dots m \dots),
P=(\dots m \dots m+1 \dots) \cr\cr
\endcases
$$

2. Column version. Assume $\s$ acts on columns.
Then
$$
G^{\s(L)}_{\theta_{P\s(Q)}}\theta_{P\s(Q)}=
\s(G^{L}_{\theta_{PQ}}\theta_{PQ})+(q^{-1}-q)E(m,L,Q)\th_0
$$
where $\th_0 \equiv \theta_{PQ}$ and
$$
E(m,L,Q)=\cases 0 & 
$ if $(L,L')=(\dots m \dots m+1 \dots),
Q=(\dots m \dots m+1 \dots) \cr
& $or $(L,L')=(\dots m+1 \dots m \dots),
Q=(\dots m+1 \dots m \dots) \cr\cr
1 & $ if  $(L,L')=(\dots m \dots m+1 \dots), 
Q=(\dots m+1 \dots  m \dots) \cr\cr
-1 & $ if $ (L,L')=(\dots m+1  \dots m \dots), 
Q=(\dots m \dots m+1 \dots) \cr\cr
\endcases
$$
}

{\bf Proof}. We will do only 2. (1. is the same). Observe that
though $G^{\s(L)}_{\theta_{P\s(Q)}}$
is a sequence of operations, it is only at most one operation that originates
the correction term $(q^{-1}-q)E(m,L,Q)\th_0$, namely 
$G_{a_{p_iq_i}a_{p_jq_j}}$, with $\{q_i,q_j\}=\{m,m+1\}$.
At this point the result is a consequence of Lemma (3.6).
The fact the correction term is not
precisely $\t_0$ but only congruent to it modulo $I_M$ comes from the
fact that the operation generating it may not be the first to occur.

\medskip

{\bf Corollary (3.12)}. {\it
Let $\s=(m,m+1) \in S_{n_0}$. 

1. Row version. Let $t_P \in T$ be a tensor in row order $P$. Then
$$
F^{\s(K)}_{t_{\s(P)}}t_{\s(P)}=
\s(F^{K}_{t_{P}}t_{P})+(q^{-1}-q)E(m,K,P)t_0
$$
where $t_0 \equiv t_{P}$ and the coefficient $E(m,K,P)$
has been defined in Lemma (3.11).

2. Column version. Let $t_Q \in T$ be a tensor in row order $Q$. Then
$$
G^{\s(L)}_{t_{\s(Q)}}t_{\s(Q)}= 
\s(G^{L}_{t_{Q}}t_{Q})+(q^{-1}-q)E(m,L,Q)t_0
$$
where $t_0 \equiv t_{Q}$ and the coefficient $E(m,L,Q)$
has been defined in Lemma (3.11).
}

{\bf Proof}. Immediate by Lemma (3.11).



\medskip

{\bf Observation (3.13)}.
In Part 1. (respectively Part 2.) one can assume to
have repeated column (respectively row) indices. The proof is the same.

\bigskip

{\bf 4. Commutation rules among the quantum minors: 
a first overview}.
\medskip 

{\bf Definition (4.1)}. Let $I=(i_1 \dots i_r)$, $1 \leq i_1 < \dots
< i_r \leq n$, $J=(j_1 \dots j_r)$, $1 \leq j_1 < \dots
< j_r \leq n$. We define $[I,J]_r \in k_q\<a_{rs}\>$ 
the {\it quantum row minor} obtained by
taking the rows $I$ and the columns $J$ as:
$$
[I,J]_r=\sum_{\s} (-q)^{-l(\s)} a_{i_1\s(i_1)}   \dots a_{i_r\s(i_r)} 
$$
where $\s$ runs over all bijective maps $\{i_1 \dots i_r\} 
\lra \{j_1 \dots j_r\}$. 

\medskip

We define $[I,J]_c \in k_q\<a_{ij}\>$
the {\it quantum column minor} obtained by
taking the rows $I$ and the columns $J$ as:
$$
[I,J]_c=\sum_{\t} (-q)^{-l(\s)} a_{\s(j_1)j_1}   \dots a_{\s(j_r)j_r}
$$
where $\s$ runs over all bijective maps $\{j_1 \dots j_r\}
\lra \{i_1 \dots i_r\}$.

\medskip

We will denote with the symbol $[I,J]$ the common image
of $[I,J]_r$ and $[I,J]_c$ in $\cO_q(M_n(k))$ under the projection
$p:k_q\<a_{ij}\> \lra  k_q\<a_{ij}\>/I_M \cong \cO_q(M_n(k))$.

We also will use the same symbol $[I,J]$ 
to denote a quantum row/column minor in $k_q\<a_{ij}\>$
whenever it is not important to distinguish
between $[I,J]_r$ or $[I,J]_c$. This happens for example when
the quantum row or column minor appears in a congruence modulo $I_M$,
since $[I,J]_r \equiv [I,J]_c$.

\medskip

For more details on quantum minors see [PW] pg 42.

\medskip

We now want to give the notion of commutation relation among quantum
minors.

\medskip

{\bf Definition (4.2)}. Let $I=(i_1 \dots i_r)$, $1 \leq i_1 < \dots
< i_r \leq n$, $J=(j_1 \dots j_r)$, $1 \leq j_1 < \dots
< j_r \leq n$, $K=(k_1 \dots k_s)$, $1 \leq k_1 < \dots
< k_s \leq n$, $L=(l_1 \dots l_s)$, $1 \leq l_1 < \dots
< l_s \leq n$. Let $(I,J) < (K,L)$, $<$ being the lexicographic
ordering.

A {\it commutation rule or relation} among 
the quantum minors $[I,J], [K,L] \in k_q\<a_{ij}\>$ 
is an expression in $k_q\<a_{ij}\>$
($1 \leq i,j \leq n$):
$$
f(q)[K,L][I,J]-g(q)[I,J][K,L]-\sum_{(I_z,J_z)<(K,L)}
b_z[I_z,J_z][K_z,L_z] \equiv 0
$$
for suitable $f(q), g(q) \in k_q$, $f(0)=g(0)=1$, $b_z\in k_q$. 
The congruence as always
has to be intended modulo $I_M$.
\medskip

We observe that since the nature of
the commutation relations among the $a_{ij}$'s (i.e. the generating
relations $R$ and $S$ of $I_M$),
the commutation relations among quantum minors will depend only
on their relative position (not the ``absolute'' position
inside the quantum matrix).

Hence without loss of generality we can assume that
given two quantum minors $[I,J]$, $[K,L]$, the largest of the
two sets $I \cup K$  and $J \cup L$  is equal to $\{ 1 \dots n\}$.

\medskip

\medskip
We now define the intertwining order of a row  multiindex $I$
and of a column multiindex $J$.
This is a central notion, since the proof of the formula for
the commutation relations is by induction on the intertwining orders.
\medskip

{\bf Definition (4.3)}. 
Let $I=(i_1 \dots i_r)$,
$1 \leq i_1 < \dots < i_r \leq n$.
Let $p$ be the first integer such that $i_p-i_{p-1} > 1$, $i_0=1$.
Define $\s_I=(i_{p}, i_{p}-1)$ , 
the {\it standard row transposition} of $I$.
Given $I$, there exists a sequence of transpositions
$\s_{I_i}$, $i=1 \dots N$, such that $I=I_N$, 
$\s_{I_i}I_i=I_{i-1}$ and $I_0=(1 \dots r)$.  
$N$ is called {\it the row intertwining order} of $I$, $intr(I)$, 
the sequence
$$
I=I_N > I_{N-1} > \dots > I_0 = (1 \dots r)
$$
is called the {\it standard row tower} of $I$ ($>$ being
the lexicographic ordering).

An example of standard row tower for $I=(34)$ is:
$$
\matrix
I=I_4=(34)>I_3=(24)>I_2=(14)>I_1=(13)>I_0=(12) \cr\cr
\s_4=(23) \qquad s_3=(12) \qquad s_2=(34) \qquad s_1=(23)
\endmatrix
$$

Let $J=(j_1 \dots j_s)$,
$1 \leq j_1 < \dots < j_s \leq n$.
Let $p$ be the first integer (starting from the
right) such that $j_p-j_{p+1} > 1$, $j_{s+1}=n$.
Define $\s_J=(j_{p}, j_{p-1}+1)$, the {\it standard transposition} of $J$.
Given $J$, there exists a sequence of transpositions
$\s_{J_i}$, $i=1 \dots N$, such that $J=J_N$, 
$\s_{J_i}J_i=J_{i-1}$ and $J_0=(n-s+1 \dots n)$.  
$N$ is called {\it the column intertwining order} of $I$, $intc(I)$, 
the sequence
$$
J=J_N < J_{N-1} < \dots < J_0=(n-s+1 \dots n)
$$
is called the {\it column standard tower} of $J$.

\medskip
An example of standard column tower for $J=(13)$ is:
$$
J=J_3=(13)<J_2=(14)<J_1=(24)<J_0=(34)
$$
We define the {\it total intertwining order $int(I,J)$} 
as the sum of the row and column intertwining orders.

\medskip

The next lemma will give the starting point for the induction.

\medskip

{\bf Lemma (4.4)}. {\bf Commutation relation, $intr=0$, $intc=0$}
{\it  
$$
[r+1 \dots n,1 \dots r][1 \dots r,r+1 \dots n] \equiv
[1 \dots r,r+1 \dots n][r+1 \dots n,1 \dots r]
$$
}

{\bf Proof.} Immediate: all elements of the first minor
commute with all the elements of the second one.

\medskip

The next lemma represents the key step in the algorithm for the
commutation relations among the quantum minors.

\medskip

{\bf Lemma (4.5)}. {\it Let $I,J,Z,W$ be ordered multiindices,
$|I|=|J|=r$, $|Z|=|W|=s$, $r+s=n$.
Let $\s=(m,m+1) \in S_n$. 

1. Row version. Assume $I \cap Z=\emptyset$
(hence $I \cup Z=\{1 \dots n\}$).
Then
$$
\align
F^{\s(K)}_{[\s(I),J]_r[\s(Z),W]_r}[\s(I),J]_r[\s(Z),W]_r &= 
\s(F^{K}_{[I,J]_r[Z,W]_r}[I,J]_r[Z,W]_r)+
\cr\cr 
&+(q^{-1}-q)E(m,K,(I,Z))t
\endalign
$$
where the coefficient $E(m,K,(I,Z))$ has been defined in
Lemma (3.11) and $t \equiv [I,J][Z,W]$.

2. Column version. Assume $J \cap W=\emptyset$ 
(hence $J \cup W=\{1 \dots n\}$). 
Then
$$
\align
G^{\s(L)}_{[I,\s(J)]_c[Z,\s(W)]_c}[I,\s(J)]_c[Z,\s(W)]_c &=
\s(G^{L}_{[I,J]_c[Z,W]_c}[I,J]_c[Z,W]_c)+
\cr\cr
&+(q^{-1}-q)E(m,L,(J,W))t
\endalign
$$
where the coefficient $E(m,L,(J,W))$ has been defined in 
Lemma (3.11) and $t \equiv [I,J][Z,W]$.
}

{\bf Proof.} This is a straightforward consequence of Corollary (3.12)
and Observation (3.13).

\medskip

{\bf Example (4.6)}. This example illustrates how to obtain
commutation relations using Lemma (4.5) in a simple
but crucial case. Consider the product of the two column minors:
$$
t=[34,13]_c[12,24]_c \in T
$$
$t$ is in column order $(1324)$. We want to apply to $t$ a sequence
of operations of type $G$, $G_t^{(24)}$, so that $G_t^{(24)}t$,
is in column order $(2413)$.
In view of an application of Lemma (4.5) we need to take:  
$(L)=(34)$, $\s(L)=(24)$, $J=(12)$, $\s(J)=(13)$, $W=(34)$, 
$\s(W)=(24)$.
$$
\align
[34,13][12,24] & \equiv 
G^{(24)}_{[34,13]_c[12,24]_c}[34,13]_c[12,24]_c  
\equiv \hbox{ by Lemma (4.5) }
\cr\cr
& \equiv \s(G^{(34)}_{[34,12]_c[12,34]_c}[34,12]_c[12,34]_c)
-(q^{-1}-q)[34,12][12,34] \equiv
\cr\cr
& \equiv  \s([12,34][34,12])
-(q^{-1}-q)[34,12][12,34] 
\cr\cr
& =[12,24][34,13]-(q^{-1}-q)[34,12][12,34]\equiv \hbox{ by Lemma (4.4) }
\cr\cr
& \equiv [12,24][34,13]-(q^{-1}-q)[12,34][34,12]
\endalign
$$
So we have obtained the following commutation relation:
$$
[34,13][12,24] \equiv [12,24][34,13]-(q^{-1}-q)[12,34][34,12] 
$$

\medskip
It is our intention to give the commutation relations for two generic
quantum minors examining all possible cases. For a complete list of
each case that arises, together with the number of
the theorem taking care
of it, see at the end of \S 7.

\bigskip
{\bf 5. Commutation relations among the quantum minors: 
the case of zero row intertwining}.
\medskip
In this section we will work out the commutation relations between
two quantum minors $[I,J]$, $[K,L]$ in the
hypothesis that $intr(I)=0$. This implies $I=(1 \dots r)$
and without loss of generality we can assume
$K=(r+1 \dots r+s)$, (with $r+s=n$). 

\medskip

We first want to take care of the situation where there are
no column indeces in common between the two quantum minors.

\medskip

Given an index $I$, let $I'$ denote its complement in 
$(1 \dots n)$. The indices are not assumed ordered if we
don't specify it.

\medskip
{\bf Theorem (5.1)}. {\it Let $I=(i_1 \dots i_r)$, 
$1 \leq i_1 <\dots <i_r \leq n$, $I'=(i'_1 \dots i'_s)$,
$1 \leq i'_1 <\dots <i'_s \leq n$,  
$n=r+s$. Let $N=intc(I)$. Then
$$
[r+1 \dots n,I'][1\dots r,I] \equiv
\sum_{i=0}^N (-1)^{i} (q^{-1}-q)^i \sum_{(R,R') \in C_i^I}
[1\dots r,R][r+1 \dots n,R']
$$
where the set $C^{I}_{i}$ is defined in the following way. 
Consider the standard column tower:
$I=I_N < \dots  < I_0=(r+1 \dots n)$,
$\sigma_N= \s =(m,m+1)$. 
Define:
$$
\matrix
C^{I_0}_{i}=((I_0,I_0')) \cr\cr 
C^{I_s}_{i}=(C^{I_{s-1}}_{i-1}-\phi^{I_{s-1}}_{i-1}, 
\sigma_s(C^{I_{s-1}}_{i})),  
\qquad 1 \leq  i \leq s \leq N 
\cr\cr 
C^{I}_{i}=C^{I_N}_{i}
\endmatrix
$$
where with the notation for $C^{I_s}_{i}$ we mean:
$$
\matrix
C^{I_s}_{i}=((X_1,X_1') \dots (X_p,X_p'),\s (Y_1,Y_1') \dots \s (Y_q,Y_q')) 
\cr \cr
C^{I_{s-1}}_{i-1}-\phi^{I_{s-1}}_{i-1}=((X_1,X_1') \dots (X_p,X_p')), \qquad
C^{I_{s-1}}_{i}=((Y_1,Y_1') \dots (Y_q,Y_q'))
\endmatrix
$$
$\phi^{I_{s}}_i$ is the subset of $C^{I_{s}}_{i}$ of
indices $(R,R')=(\dots z \dots z+1 \dots )$, $\s_s=(z,z+1)$.

We agree to put $C^{I_{s}}_{i}$ and
$\phi^{I_{s}}_i$ equal to the empty set if one of the indices $i$ or $s$
is negative or if $i > s$. 
Notice that the $(X_j,X_j')$'s and $\s (Y_l,Y_l')$'s 
in $C_i^{I_s}$ need not to be distinct.
}

{\bf Proof.} This proof is based on the same idea as the proof
of Theorem (2.19) in [Fi1] pg 429. We however 
will write the essential steps since here the context is
different and some sign appears.

The proof is by induction on $N$ the column intertwining order of $I$.
When $N=0$ this is Lemma (4.4).

\medskip

Let $G^I_{[r+1\dots n,I']_c[1\dots r,I]_c}$ 
be the sequence of operations of type $G$ such that  \break
$G^I_{[r+1\dots n,I']_c[1\dots r,I]_c}[r+1\dots n,I']_c[1\dots r,I]_c$ 
is in column order $(I,I')$.  
This will lead to the commutation of the minors.

\medskip
Lemma (4.5) allows us to write:
$$
\matrix
[r+1\dots n,I'][1\dots r,I]  \equiv
G^I_{[r+1\dots n,I']_c[1\dots r,I]_c}[r+1\dots n,I']_c[1\dots r,I]_c=  
\cr\cr
=\s(G^{I_{N-1}}_{[r+1\dots n,I_{N-1}']_c[1\dots r,I_{N-1}]_c}
[r+1\dots n,I_{N-1}']_c[1\dots r,I_{N-1}]_c)
-(q^{-1}-q)t_{I_{N-1}}
\endmatrix
$$
where by Lemma (4.5)
$$
\align
t_{I_{N-1}} & \equiv [r+1\dots n,I'_{N-1}][1\dots r,I_{N-1}]  \equiv 
\hbox{ by induction }
\cr\cr
& \equiv
\sum_{i=0}^{N-1} (-1)^{i}(q^{-1}-q)^i \sum_{(R,R') \in C_i^{I_{N-1}}}
[1\dots r,R][r+1 \dots n,R']
\endalign
$$
By Lemma (3.11), the coefficient $E(m,I_{N-1}, (I_{N-1}',I_{N-1}))=-1$ since
$(L,L')=(I_{N-1},I_{N-1}')=( \dots m+1 \dots m \dots)$, while
$(J,W)=(I_{N-1}',I_{N-1})=( \dots m \dots m+1 \dots)$, because
a standard column tower has descending lexicographic ordering
i.e. $I_N<I_{N-1}$ and $\s(I_{N-1})=I_N$, $\s=(m,m+1)$.

\medskip

To finish our proof we just need to show that:
$$
\matrix
\s(G^{I_{N-1}}_{[r+1\dots n,I'_{N-1}][1\dots r,I_{N-1}]}
[r+1\dots n,I'_{N-1}][1\dots r,I_{N-1}])
\equiv 
\cr\cr
\equiv
\sum_{i=0}^{N-1} (-1)^{i}(q^{-1}-q)^i \sum_{(R,R') \in \s(C_i^{I_{N-1}})}
[1\dots r,R][r+1 \dots n,R']
\cr\cr
-\sum_{i=0}^{N}(-1)^{i}(q^{-1}-q)^i\sum_{(R,R') \in \phi_{i-1}^{I_{N-1}}}
[1\dots r,R][r+1 \dots n,R']
\endmatrix
$$

\medskip
Let
$$
t_0=\sum_{i=0}^{N-1}(-1)^{i}(q^{-1}-q)^i\sum_{(R,R') \in C_{i}^{I_{N-1}}}
[1\dots r,R][r+1 \dots n,R']
$$
By induction hypothesis we have that
$$
G^{I_{N-1}}_{[r+1\dots n,I_{N-1}'][1\dots r,I_{N-1}]}
[r+1\dots n,I_{N-1}'][1\dots r,I_{N-1}]
\equiv t_0\equiv G^{I_{N-1}}_{t_0}t_0
$$ 
But actually we have
$$
G^{I_{N-1}}_{[r+1\dots n,I_{N-1}'][1\dots r,I_{N-1}]}
[r+1\dots n,I_{N-1}'][1\dots r,I_{N-1}]=G^{I_{N-1}}_{t_0}t_0
$$
by Lemma (2.12) in [Fi1].

So we are left with the easier problem to compute
$\s(G^{I_{N-1}}_{t_0}t_0)$.

\medskip

Write $t_0=t_1+t_2$, where
$$
\matrix
t_1=\sum_{i=0}^{N-1}(-1)^{i}(q^{-1}-q)^i\sum_{(R,R') \in C_{i}^{I_{N-1}}
-\phi^{I_{N-1}}_i}
[1\dots r,R][r+1 \dots n,R'] \cr\cr
t_2=\sum_{i=0}^{N-1}(-1)^{i}(q^{-1}-q)^i\sum_{(R,R') \in \phi_{i}^{I_{N-1}}}
[1\dots r,R][r+1 \dots n,R']
\endmatrix
$$

Let's put $\s(t_0)$ in order $(I,I')$. Using Lemma (4.5) one of
the terms will come out $\s(G^{I_{N-1}}_{t_0}t_0)$.
$$
\align
\s(t_0)&=\s(t_1)+\s(t_2) \equiv G^I_{\s(t_1)}\s(t_1)+
G^I_{\s(t_2)}\s(t_2) =
\cr\cr
&=G^{\s(I_{N-1})}_{\s(t_1)}\s(t_1)+G^{\s(I_{N-1}) }_{\s(t_2)}\s(t_2)\equiv 
\hbox{ by Corollary (3.12) }
\cr\cr
&\equiv\s(G^{I_{N-1}}_{t_1}t_1)+\s(G^{I_{N-1}}_{t_2}t_2) 
\cr\cr
&-(q^{-1}-q)\sum_{i=0}^{N-1}(-1)^{i}(q^{-1}-q)^i
\sum_{(R,R') \in \phi_{i}^{I_{N-1}}}
[1\dots r,R][r+1 \dots n,R'] =
\cr\cr
&=\s(G^{I_{N-1}}_{t_0}t_0)+
\cr\cr
&+\sum_{i=0}^{N-1}(-1)^{i+1}(q^{-1}-q)^{i+1}
\sum_{(R,R') \in \phi_{i}^{I_{N-1}}}
[1\dots r,R][r+1 \dots n,R']=
\cr\cr
&=\s(G^{I_{N-1}}_{t_0}t_0)-
\sum_{i=1}^{N}(-1)^{i}(q^{-1}-q)^i
\sum_{(R,R') \in \phi_{i-1}^{I_{N-1}}}
[1\dots r,R][r+1 \dots n,R'] 
\endalign
$$
By Lemma (3.11), the coefficient $E(m,I_{N-1}, (R,R'))=-1$
for $(R,R') \in \phi_i^{I_{N-1}}$. In fact 
$(I_{N-1},I_{N-1}')=( \dots m+1 \dots m \dots )$ and
$\phi_i^{I_{N-1}}$ is precisely defined as the subset of
$C_i^{I_{N-1}}$ containing those indices $(R,R')=( \dots m \dots m+1 \dots)$.

Consequently the coefficient $E(m,I_{N-1}, (R,R'))=0$ 
for $(R,R') \in C_i^{I_{N-1}} \-\phi_i^{I_{N-1}}$.

\medskip

Hence we have:
$$
\align
\s(G^{I_{N-1}}_{t_0}t_0) &\equiv\s(t_0)+
\sum_{i=1}^{N}(-1)^{i}(q^{-1}-q)^i
\sum_{(R,R') \in \phi_{i-1}^{I_{N-1}}}
[1\dots r,R][r+1 \dots n,R']
\cr\cr
&\equiv 
\s(\sum_{i=0}^{N-1}(-1)^{i}(q^{-1}-q)^i\sum_{(R,R') \in C_{i}^{I_{N-1}}}
[1\dots r,R][r+1 \dots n,R'])+
\cr\cr
&-\sum_{i=0}^{N}(-1)^{i}(q^{-1}-q)^i
\sum_{(R,R') \in \phi_{i-1}^{I_{N-1}}}
[1\dots r,R][r+1 \dots n,R']
\endalign
$$
which is what we wanted to show.

\medskip
Given a generic multiindex $S=(s_1 \dots s_r)$, $1 \leq s_1 \dots s_r \leq n$,
we define $S_{ord}=(s_1^{ord} \dots s_r^{ord})$,
$\{s_1 \dots s_r\}=\{s_1^{ord} \dots s_r^{ord}\}$ and 
$1 \leq s_1^{ord} < \dots < s_r^{ord} \leq n$. 
We also define $l(S)$ the
length of the permutation sending $S$ into $S_{ord}$.

\medskip
{\bf Corollary (5.2): Commutation relations for $intr$=0 and no repeated 
indices}. 

{\it Let $I$ and its complement $I'$ in $\{1 \dots n\}$ be two
ordered indices. Then:
$$
\matrix
[r+1 \dots n,I'][1\dots r,I] \equiv
[1\dots r,I][r+1 \dots n,I']+
\cr\cr
+\sum_{i=1}^N (-1)^{i} (q^{-1}-q)^i \sum_{(R,R') \in C_i^I}
(-q)^{-l(R)-l(R')}[1\dots r,R][r+1 \dots n,R']
\endmatrix
$$
is a commutation relation between the two
quantum minors $[r+1 \dots n, I']$ and $[1 \dots r,I]$.
}

{\bf Proof}. Immediate from (5.1) and the fact that
$[I,J] \equiv (-q)^{-l(J)}[I,J_{ord}]$ (see [PW] pg 42). 

\medskip

{\bf Example (5.3)}. 
In $\cO_q(M_4(k))$ let's consider the commutation relation between
$[34,34]$, $[12,12]$. According to the theorem we have to build
a standard column tower:
$$
\matrix
I_4=(12) > I_3=(13) > 
I_2=(14) > I_1=(24) > I_0=(34)  \cr\cr
\s_4=(23), \qquad 
\s_3=(34), \qquad
\s_2=(12), \qquad \s_1=(23).
\endmatrix
$$
We now build the sets $C_i^{I_s}$. To improve readability we
will write not only the column index, but also the row one.
So for example instead of writing under the column of
$C_1^{I_2}$, $(24,13)$, $(34,21)$ we will write:
$(12,24)(34,13)$, $(12,34)(34,21)$.

\medskip

$$
\matrix
s & \s_s& C_0^{I_s} & C_1^{I_s} & C_2^{I_s} & C_3^{I_s} \cr\cr
0 & & (12,34)(34,12) & & & \cr\cr\cr
1 & (23)& (12,24)(34,13) & (12,34)(34,12) & && \cr\cr
2 & (12)& (12,14)(34,23) & (12,24)(34,13) & (12,34)(34,12) \NO & \cr\cr
	&&&			(12,34)(34,21) & &
\cr\cr\cr
3 & (34)& (12,13)(34,24) & (12,23)(34,14) & (12,24)(34,13)   & \cr\cr
&&&				(12,14)(34,23) & (12,43)(34,12) \NO & \cr\cr
&&&				(12,43)(34,21) && 
\cr\cr\cr
4 & (23)& (12,12)(34,34) & (12,13)(34,24) & (12,23)(34,14)\NO & 
(12,24)(34,13)\NO\cr\cr
&&			      &	(12,32)(34,14) & (12,14)(34,23)\NO & \cr\cr
&&			      &	(12,14)(34,32) & (12,43)(34,21) & \cr\cr
&&			      &	(12,42)(34,31) & (12,34)(34,12) & \cr\cr
\endmatrix
$$
The commutation relation is:
$$
\align
[34,34][12,12] & \equiv[ 12,12][34,34]+
\cr\cr
&-(q^{-1}-q)[[12,13][34,24]+
[12,32][34,14]+[12,14][34,32]+
\cr\cr
&+[12,42][34,31]]+\cr\cr
&+(q^{-1}-q)^2[[12,43][34,21]+[12,34][34,12]]
\endalign
$$
which becomes:
$$
\align
[34,34][12,12]& \equiv [12,12][34,34]+
\cr\cr
&-(q^{-1}-q)[[12,13],[34,24]-q^{-1}[12,23][34,14]+
\cr\cr
&-q^{-1}[12,14][34,23]
+q^{-2}[12,24][34,13]]
\cr\cr
&+(q^{-1}-q)^2(1+q^{-2})[12,34][34,12]
\endalign
$$
as one can directly check.

\medskip
We now want to consider the case of repeated column indeces.

\medskip
{\bf Proposition (5.4)}. {\it Let $J=(j_1 \dots j_r)$ and
$L=(l_1 \dots l_s)$ be two ordered indices, $r+s=n$.
Let $[1 \dots r,J]$, 
$[r+1 \dots n,L]$ be two quantum
minors, with $C=(c_1 \dots c_z)$ common column indices.
Assume $l_u<j_v$ for all $l_u \in L-C$
and $j_v \in J-C$.
Then
$$
q^{-z}[r+1 \dots n,l_1 \dots l_s]
[1 \dots r,j_1 \dots j_r] \equiv
[1 \dots r,j_1 \dots j_r]
[r+1 \dots n,l_1 \dots l_s]
$$ 
}

{\bf Proof.} Direct calculation based on [Fi1] Lemma (2.11).

\medskip

{\bf Theorem (5.5). Commutation relation $intr$=0
(no repeated row indices).}

{\it Let $I=(i_1 \dots i_r)$, 
$1 \leq i_1 < \dots < i_r \leq n$, $J=(j_1 \dots j_s)$,
$1 \leq j_1 < \dots < j_s \leq n$, $I \cap J= C=(c_1 \dots c_z)$, 
$I \cup J=\{1 \dots n\}$. 
Let $N=intc(I \- C)$ the length of the column standard
tower of $I \- C$ in the set $(I\cup J) \- C$. Then
$$
\matrix
q^{-z}[r+1 \dots n,J][1\dots r,I] \equiv
\sum_{i=0}^N (-q)^i(q^{-1}-q)^i \cr\cr
\sum_{(R,R') \in C_i^{I \- C}} (-q)^{-l(R)-l(R')}
[1\dots r,(R \cup C)_{ord}][r+1 \dots n,(R' \cup C)_{ord}]
\endmatrix
$$
where the complement $R$ of an column index $R'$ belonging to 
the standard tower of $I\- C$ is taken in $I \cup J \- C$.
}

{\bf Proof.} 
Let's write the quantum minors using the Laplace expansion as
in [PW] pg 46. According to the notation in [PW] a
multiindex $X=(x_1 \dots x_p)$ is also seen as an integer 
$x_1 + \dots + x_p$. 
$$
\matrix
q^{-z}[\{r+1 \dots n\},J][\{1 \dots r\},I] \equiv
q^{-z}\sum_{S} q^{C-S}[S,C][\{r+1 \dots n\} \- S,J \- C]
\cr\cr
\sum_{\tS} q^{\tS-C}[\{1 \dots r\} \- \tS,I \- C][\tS,C]=
\cr\cr
=q^{-z}\sum_{S,\tS}q^{\tS-S}[S,C][\{r+1 \dots n\} \- S,J \- C]
[\{1 \dots r\} \- \tS,I \- C][\tS,C]
\endmatrix
$$
where the sums are over any $S \subset \{r+1 \dots n\}$ with  
$z=|C|$ elements and over any $\tS \subset \{1 \dots r\}$
with $z=|C|$ elements. 

Since $(J \- C) \cap  (I \- C) = \emptyset$,
by Theorem (5.1) we can write
$$
\matrix
[\{r+1 \dots n\} \- S,J \- C]
[\{1 \dots r\} \- \tS, I \- C] \equiv \cr\cr
\equiv \sum_{i \in [0,N], (R,R') \in C_i^{I \- C}} \a_i 
[\{1 \dots r\} \- \tS,R][\{r+1 \dots n\}\- S,R']
\endmatrix
$$
with $\a_i=(-q)^i(q^{-1}-q)^i$, $R'$ is the complement of
the element $R$ in the standard tower of $I \- C$ in
$(I \cup J )\- C$.

Let's substitute in the previous expression:
$$
\matrix
q^{-z}[\{r+1 \dots n\},J][\{1\dots r\},I]  \equiv
q^{-z}\sum_{S,\tS}
\sum_{i \in [0,N], (R,R') \in C_i^{I \- C}} \a_i
q^{\tS-S}[S,C] \cr\cr
[\{1 \dots r\} \- \tS,R][\{r+1 \dots n\}\- S,R']
[\tS,C]
\endmatrix
$$
We now need to commute $[S,C]$ with $[\{1 \dots r\} \- \tS,R]$
and $[\{r+1 \dots n\}\- S,R']$ with $[\tS,C]$. Observe that
$C \cap R = \emptyset$ and $C \cap R' = \emptyset$, hence
we can again apply Theorem (5.1). 
For the first commutation rule we need to build a standard
tower of $R$ in $C \cup R$.
$$
\matrix
[S,C][\{1 \dots r\} \- \tS,R]  \equiv 
\sum_{i=0}^N \sum_{(U,U') \in C_i^{R}} \a_i 
[\{1 \dots r\} \- \tS,U][S,U']
\endmatrix
$$
Similarly for the second commutation rule, we need to build a standard
tower of $C$ in $R' \cup C$:
$$
\matrix
[\{r+1 \dots n\}\- S,R'][\tS,C]  \equiv 
\sum_{i=0}^N\sum_{ (V,V') \in C_i^{C}} \a_i
[\tS,V][\{r+1 \dots n\}\- S,V']
\endmatrix
$$
Let's substitute the first expression for the commutation rule.
$$
\matrix
q^{-z}[\{r+1 \dots n\},J][\{1 \dots r\},I]  \equiv
q^{-z}\sum_{S,\tS}
\sum_{i \in [0,N], (R,R') \in C_i^{I \- C}} \a_i
q^{\tS-S}\cr\cr
\sum_{j=0}^N \a_i
[\{1 \dots r\} \- \tS,R][S,C]
[\{r+1 \dots n\}\- S,R']
[\tS,C]+
\cr\cr
+q^{-z}\sum_{S,\tS}
\sum_{i \in [0,N], (R,R') \in C_i^{I \- C}} \a_i
q^{\tS-S}\cr\cr
\sum_{j=1}^N \sum_{(U,U') \in C_j^{R}} \a_j
[\{1 \dots r\} \- \tS,U][S,U']
[\{r+1 \dots n\}\- S,R']
[\tS,C]
\endmatrix
$$
Observe that by the very construction $U'\neq R'$.
Hence the sum:
$$
\sum_{S}(-q)^{R-S} [S,U']
[\{r+1 \dots n\}\- S,R'] \equiv 0
$$
by [PW] pg 46. So we are left with:
$$
\align
q^{-z}[\{r+1 \dots n\},J][\{1 \dots r\},I]  \equiv &
q^{-z}\sum_{S,\tS}
\sum_{i \in [0,N], (R,R') \in C_i^{I \- C}} \a_i
q^{\tS-S}\cr\cr
& [\{1 \dots r\} \- \tS,R][S,C]
[\{r+1 \dots n\}\- S,R']
[\tS,C]
\endalign
$$
Reasoning in the same way for the commutation of
$[\{r+1 \dots n\}\- S,R']
[\tS,C]$ we get by Proposition (5.4):
$$
\matrix
q^{-z}[\{r+1 \dots n\},J][\{1 \dots r\},I]  \equiv
q^{-z}\sum_{S,\tS}
\sum_{i \in [0,N], (R,R') \in C_i^{I \- C}} \a_i
q^{\tS-S} \cr\cr
[\{1 \dots r\} \- \tS,R][S,C]
[\tS,C][\{r+1 \dots n\}\- S,R'] \equiv
\cr\cr
\equiv
\sum_{S,\tS}
\sum_{i \in [0,N], (R,R') \in C_i^{I \- C}} \a_i
q^{\tS-S}
[\{1 \dots r\} \- \tS,R][\tS,C][S,C][\{r+1 \dots n\}\- S,R']
\endmatrix
$$
which is after summing over $S$ and $\tilde S$ yields our result.

\medskip

{\bf Example (5.6)}. Consider the two quantum minors
$[34,23]$, $[12,13]$. $C=(3)$. A standard column tower for
$I \- C = (2)$ is $(1) < (2)$, $\s =(12)$. 
$I \cup J \- C=\{1,2\}$.

The sets $C_i^{I\-C}$ are:
$$
\matrix
s & \s_s& C_0^{(I-C)_s} & C_1^{(I-C)_s}  \cr\cr
0 & & (12,2)(34,1) &  \cr\cr
1 & (12) & (12,1)(34,2) & (12,2)(34,1)
\endmatrix
$$
Hence the commutation is
$$
q^{-1}[34,23][12,13] \equiv
[12,13][34,23]-(q^{-1}-q)[12,23][34,13]
$$
as one can directly check.

\medskip



\bigskip

{\bf 6. Commutation rules among the quantum minors: the case of
no row indices in common}.
\medskip 
Assume first that given two
quantum minors $[I,J]$, $[K,L]$ 
$I \cap K= \emptyset$, $J \cap L=\emptyset$, i.e. the
two minors have no row or column indices in common.  So we
can assume without loss of generality that 
$I \cup K=\{1 \dots n\}$ ($J \cup L=\{1 \dots n\}$). 

\medskip
Given 
and multiindex $I$, let $I'$ denote as always its complement
in $\{1 \dots n\}$.

\medskip

{\bf Theorem (6.1)}. {\it Let $I,J$ be ordered multiindices,
$|I|=r$, $|J|=s$, $r+s=n$.
Then
$$
[I',J'][I,J]
\equiv
\sum_{i=0}^{M+N} (q^{-1}-q)^i \sum_{(Z,W) \in C_i^{I,J}}
sign(Z,W) [Z,W][Z',W']
$$
where the set $C^{I,J}_{i}$ is defined in the following way. 
Consider the standard column tower:
$$
J=J_M < \dots  < J_0=(r+1 \dots n)
$$
and the standard row tower
$$
I=I_N > \dots > I_0=(1 \dots r)
$$ 
Let $\s_s^c$, $\s_s^r$ denote respectively the permutations
relative to the $s$ step of the column, respectively row, tower.

Define:
$$
\matrix
C^{I_0,J_0}_{0}=((I_0,J_0)) \cr 
C^{I_0,J_s}_{i}=(C^{I_0,J_{s-1}}_{i-1}-\phi^{I_0,J_{s-1}}_{i-1}, 
\sigma_s^c(C^{I_0,J_{s-1}}_{i})),  
\qquad 1 \leq  i \leq s \leq M 
\cr \cr
C^{I_t,J_M}_{i}=(C^{I_{t-1},J_{M}}_{i-1}-\psi^{I_{t-1},J_{M}}_{i-1},
\sigma_t^r(C^{I_{t-1},J_{M}}_{i})),
\qquad 1 \leq  i \leq t \leq N 
\cr\cr
C^{I,J}_{i}=C^{I_N,J_M}_{i}
\endmatrix
$$
$\phi^{I_{0},J_s}_i$ is the subset of $C^{I_0,J_{s}}_{i}$ of
indices $(I_0,U)$ such that
$(U,U')=(\dots m \dots m+1 \dots)$, with $\s_s^c=(m,m+1)$.

$\psi^{I_{t},J_M}_i$ is the subset of $C^{I_t,J_{M}}_{i}$ of
indices $(V,V')=(\dots m+1 \dots m \dots )$ such that, $\s_t^r=(m,m+1)$.

The sign of $(Z,W)$ is defined as follows.
$$
sign(Z,W)=(-1)^i \qquad \hbox{ for all }  
(Z,Z') \in C^{I_0,J_s}_i.
$$
For $(Z,W) \in C^{I_t,J_M}_i$, if $(Z,W) \in C^{I_{t-1},J_M}_{i-1}$
the sign is already defined. 
If $(Z,W) =\s^r_t((\tilde Z, \tilde W)) 
\in \s^r_t(C^{I_{t-1},J_M}_{i})$ we take 
$sign(Z,W)=sign(\tilde Z, \tilde W)$.

We take $C_i^{I_t,J_s}=0$ if $i$ or $t$ or $s$ are negative 
or if $i>s$, $i>t$.
}

{\bf Proof}. This proof is based on the same idea as proof of Theorem (5.1). 

By induction on the row intertwining order, $intr(I)$.
If $intr(I)=0$ this is Theorem (5.1).  
Let $N=intr(I)>0$. Let $\s=\s_N^r=(m,m+1)$.

By Lemma (4.5) we have:
$$
\matrix
[I',J'][I,J]  \equiv
F^I_{[I',J']_r[I,J]_r}[I',J']_r[I,J]_r=
\cr\cr
=\s(F^{I_{N-1}}_{[I_{N-1}',J']_r[I_{N-1},J]_r}
[I_{N-1}',J']_r[I_{N-1},J]_r)
+(q^{-1}-q)t_{I_{N-1}}
\endmatrix
$$
The coefficient $E(m,I_{N-1}, (I_{N-1}',I_{N-1}))=1$ since
$(I_{N-1},I_{N-1}')=( \dots m \dots m+1 \dots)$, while
$(I,Z)=(I_{N-1}',I_{N-1})=( \dots m+1 \dots m \dots)$, because
a standard row tower has descending lexicographic ordering.

By Lemma (4.5) and by our induction hypothesis
$$
\align
t_{I_{N-1}} & \equiv [I',J'][I,J] \equiv
\cr\cr 
& \equiv
\sum_{i=0}^{M+N-1} (q^{-1}-q)^i \sum_{(Z,W) \in C_i^{I_{N-1},J}}
sign(Z,W) [Z,W][Z',W']=:t_0
\endalign
$$

To finish the proof it is enough to show that
$$
\matrix
\s(F^{I_{N-1}}_{[I_{N-1}',J']_r[I_{N-1},J]_r}
[I_{N-1}',J']_r[I_{N-1},J]_r)=
\cr\cr
\sum_{i=0}^{M+N-1} (q^{-1}-q)^i \sum_{(Z,W) \in C_i^{I_{N-1},J}}
sign(Z,W) [Z,W][Z',W']+
\cr\cr
-\sum_{i=0}^{M+N} (q^{-1}-q)^i \sum_{(Z,W) \in \psi_{i-1}^{I_{N-1},J}}
sign(Z,W) [Z,W][Z',W']
\endmatrix
$$
Reasoning as in the proof of Theorem (5.1) one can see that
$$
F^{I_{N-1}}_{[I_{N-1}',J']_r[I_{N-1},J]_r}
[I_{N-1}',J']_r[I_{N-1},J]_r=
F^{I_{N-1}}_{t_0}t_0
$$
So we have to compute $\s(F^{I_{N-1}}_{t_0}t_0)$.
$$
\align
\s(t_0) & \equiv F^I_{\s(t_0)}\s(t_0) \equiv \hbox{ by Lemma (4.5) } \cr\cr
&\equiv \s(F^{I_{N-1}}_{t_0}t_0)+(q^{-1}-q)
\sum_{i=0}^{M+N-1} (q^{-1}-q)^i  \cr\cr
& \sum_{(Z,W) \in \psi_{i}^{I_{N-1},J}}
sign(Z,W) [Z,W][Z',W']
\endalign
$$
The coefficient $E(m,I_{N-1}, (Z,Z'))=+1$
for $(Z,Z') \in \psi_i^{I_{N-1}}$. In fact \break 
$(I_{N-1},I_{N-1}')=( \dots m \dots m+1 \dots )$ and
$\psi_i^{I_{N-1},J}$ is precisely defined as the subset of
$C_i^{I_{N-1},J}$ containing those indices 
$(Z,W)$, such that $(Z,Z')=( \dots m+1 \dots m \dots)$.

So we have:
$$
\align
\s(F^I_{\s(t_0)}\s(t_0))  \equiv &
\s(t_0)+ \cr\cr
& -\sum_{i=0}^{M+N-1} (q^{-1}-q)^{i+1} \sum_{(Z,W) \in \psi_{i}^{I_{N-1},J}}
sign(Z,W) [Z,W][Z',W'] =
\cr\cr
=&  \s(\sum_{i=0}^{M+N-1} (q^{-1}-q)^i \sum_{(Z,W) \in C_i^{I_{N-1},J}}
sign(Z,W) [Z,W][Z',W'])+
\cr\cr
&-\sum_{i=1}^{M+N} (q^{-1}-q)^{i} \sum_{(Z,W) \in \psi_{i-1}^{I_{N-1},J}}
sign(Z,W) [Z,W][Z',W'] 
\endalign
$$
which is what we wanted to prove.

\medskip
{\bf Example (6.2)}
Consider the commutation of the quantum minors:
$[23,13]$, $[14,24]$.
The standard row and columns
towers are always built considering the second minor.
To each transposition we will add a suffix $r$ or $c$ depending if
it belongs to a row or a column tower, i.e. it acts on row or
column indices.

\medskip
Standard row tower:
$$
(14)>_{(34)^r} (13)>_{(23)^r} (12)
$$

Standard column tower:
$$
(24)<_{(23)^c} (34)
$$

Let's compute the sets $C_i^{I_t,J_s}$. Instead of just writing 
$(Z,W)$ for an element of $C_i^{I_t,J_s}$, we will write $(Z,W)(Z',W')$.
In front of each element of $C_i^{I_t,J_s}$ we will put its sign.
$$
\matrix
\s_s& C_0^{I_t,J_s} & C_1^{I_t,J_s} & C_2^{I_t,J_s} & C_3^{I_t,J_s} 
\cr\cr
& +(12,34)(34,12) & & & \cr\cr\cr
(23)^c& +(12,24)(34,13) & -(12,34)(34,12) & & \cr\cr
(23)^r& +(13,24)(24,13) & +(12,24)(34,13) & -(12,34)(34,12) & \cr\cr
	&&		     -(13,34)(24,12) & &
\cr\cr\cr
(34)^r& +(14,24)(23,13) & +(13,24)(24,13) & +(12,24)(34,13)   & 
-(12,34)(34,12)
\cr\cr
&&			     +(12,24)(43,13) & -(13,34)(24,12) &   
\cr\cr
&&			     -(14,34)(23,12) & -(12,34)(43,12) &
\cr\cr\cr
\endmatrix
$$
The relation is:
$$
\align
[23,13][14,24]& \equiv  [14,24][23,13] +\cr\cr
&+(q^{-1}-q)[[13,24][24,13]
-q^{-1}[12,24][34,13]-[14,34][23,12]]\cr\cr
&+(q^{-1}-q)^2[[12,24][34,13]-[13,34][24,12]+q^{-1}[12,34][34,12]]+
\cr\cr
&-(q^{-1}-q)^3[12,34][34,12]
\endalign
$$
as one can check directly.
\medskip
We now want to consider the commutation of two quantum minors
in the case of positive row and column intertwining and repeated 
column indices.
\medskip
{\bf Theorem (6.3)} {\bf Commutation relation, no repeated row indices}.
{\it Let $I$, $J$, $K$, $L$ be ordered indices,
$|I|=|K|=r$, $|J|=|L|=s$, $r+s=n$.
$I \cap K= \emptyset$, $J \cap L=C$, $|C|=z$. 
Then
$$
\align
q^{-z}[K,L][I,J]
\equiv &
\sum_{i=0}^{M+N} (q^{-1}-q)^i \sum_{(Z,W) \in C_i^{I,J \- C}}
sign(Z,W) (-q)^{-l(W)-l(W')} \cr\cr
&[Z,(W \cup C)_{ord}][Z',(W' \cup C)_{ord}]
\endalign
$$
where the set $C^{I,J \- C}_{i}$ have been defined above.
$N$ is the lenght of the standard row tower of $I$ in $\{1 \dots n\}$,
$M$ is the length of the standard column tower of $J \- C$ in
$J \cup L \- C$. $W'$ is the complement of $W$ in $J \cup L \- C$.
}

{\bf Proof.} By induction on $N=intr(I)$. If $N=0$ this is Theorem (5.5).
For positive row intertwining, the argument is the same as Theorem (6.1).
We want to remark that though Lemma (4.5) (on which the proof of Theorem (6.1)
is based) is stated for tensor monomials with no repeated indices, we
noticed in Observation (3.13) that this assumption is not essential.

\bigskip

{\bf 7. Commutation rules among quantum minors: the remaining cases.}
\medskip

In this section we want to take care of the remaining cases for the
commutation relations among quantum minors to cover
all the possibilities that can arise.
\medskip
For completeness we want to give the formulas also in these
cases, though some descend almost
immediately from the previous ones.
\medskip

Let's first examine the case of repeated row indices and no repeated
column indices. 
\medskip

{\bf Theorem (7.1). Commutation relation
$intc$=0 (no repeated column indices).}

{\it Let $I=(i_1 \dots i_r)$, 
$1 \leq i_1 < \dots < i_r \leq n$, $J=(j_1 \dots j_s)$,
$1 \leq j_1 < \dots < j_s \leq n$, $I \cap J= R=(r_1 \dots r_z)$, 
$\{1 \dots n\}=\{1 \dots r\} \cup \{1 \dots s\}$. 
Let $N=intc(I \- R)$ the length of the row standard
tower of $I \- R$ in the set $(I\cup J) \- R$. Then
$$
\align
q^{-z}[J,r+1 \dots n][I,1\dots r] \equiv&
\sum_{i=0}^N (q^{-1}-q)^i 
\sum_{(Z,Z') \in C_i^{I \- R}} (-q)^{-l(Z)-l(Z')}
\cr\cr
&[(Z \cup R)_{ord},1\dots r][(Z' \cup R)_{ord},r+1 \dots n]
\endalign
$$
where the complement $Z$ of an index $Z'$ belonging to 
the standard tower of $I \- R$ is taken in $(I \cup J) \- R$.
}

{\bf Proof}. Same argument as for the column version.

\medskip
{\bf Theorem (7.2). Commutation relation
with no repeated column indices.}
{\it Let $I$, $J$, $K$, $L$ be ordered indices,
$|I|=|K|=r$, $|J|=|L|=s$, $r+s=n$.
$I \cap K= R$, $J \cap L=\emptyset$, $|R|=z$. 
Then
$$
\align
q^{-z}[K,L][I,J]
\equiv &
\sum_{i=0}^{M+N} (q^{-1}-q)^i \sum_{(Z,W) \in C_i^{I \- R,J}}
sign(Z,W) (-q)^{-l(Z)-l(Z')} 
\cr\cr
& [(Z \cup R)_{ord},W][(Z' \cup R)_{ord},W']
\endalign
$$
where the set $C^{I \- R,J}_{i}$ have been defined above.
$M$ is the lenght of the standard column tower of $J$ in $\{1 \dots n\}$,
$N$ is the length of the standard row tower of $I \- R$ in
$(I \cup K)\-  R$. $Z'$ is the complement of $Z$ in $(I \cup K) \- R$.
}

{\bf Proof}. Same arguments as for the column version (Theorem (5.5)).

\medskip

We now want to examine 
the case of two quantum minors $[I,J]$, $[K,L]$ with repeated row
and column indices, i.e. $I \cap K \neq \emptyset$,
$J \cap L \neq \emptyset$. This implies that the two quantum matrices
whose minors we are computing have at least one element in common.

\medskip
Let $S$ denote the antiendomorphism of $\cO_q(M_n(k))$ ([PW] pg 55):
$S(x_{ij})=(-q)^{j-i}[\{1 \dots n\} \- j, \{1 \dots n\} \- i]$. By Lemma (3.1)
in [Fi1] we have that: 
$$
S([I,J])=(-q)^{J-I}[\{1 \dots n\} \- I, \{1 \dots n \}\- J]
$$
\medskip

To compute a commutation relation among $[I,J]$ and $[K,L]$ 
it is now enough to apply $S$ the resulting expression will
have disjoint row and/or column indices and then we
can apply one of the theorems above. 

\medskip

We now assume to be in $\cO_q(M_n(k))$.

\medskip

{\bf Theorem (7.3)}. {\bf Commutation
of quantum minors with repeated row and column indices}
{\it Let $I$, $J$, $K$, $L$, be ordered
indices.  Let $I \cap K=R \neq \emptyset$,
$|R|=\rho$, 
$J \cap L =C \neq \emptyset$, $|C|=\gamma$. 

i) If $|I \cup K| \geq |J \cup L|$
then
$$
\align
q^{-z}[K,L][I,J]=&
\sum_{i=0}^{M+N} (q^{-1}-q)^i \sum_{(Z,W) \in C_i^{I \- R,J \- C}}
sign(Z,W) (-q)^{-l(W)-l(W')} \cr\cr
&[(Z \cup R)_{ord},(W \cup C)_{ord}][(Z'\cup R)_{ord},(W' \cup C)_{ord}]
\endalign
$$
where $M$ and $N$ are the lenghts of standard row and column  towers
of $I \- R$ and $J \- C$ in respectively $(I \cup K) \- R$ and 
$(J \cup L) \- C$.
$z=\gamma - \rho$.

ii) If $|J \cap L| \geq |I \cap K|$
then
$$
\align
q^{-z}[K,L][I,J]=&
\sum_{i=0}^{M+N} (q^{-1}-q)^i \sum_{(Z,W) \in C_i^{I \- R,J \- C}}
sign(Z,W) (-q)^{-l(Z)-l(Z')}
\cr\cr
& [(Z \cup R)_{ord},(W \cup C)_{ord}][(Z' \cup R)_{ord},(W' \cup C)_{ord}]
\endalign
$$
where $M$ and $N$ are the lenghts of standard row and column  towers
of $I \- R$ and $J \- C$ in respectively $(I \cup K) \- R$ and 
$(J \cup L) \- C$.
$z=\rho-\gamma$.  
}

{\bf Proof}. Let's briefly sketch $(i)$ ($(ii)$ being the same).

Consider the antiendomorphism $S$ in
matrix bialgebra generated by the elements $x_{ij}$
with $i \in I \cup K$, $j \in J \cup L$ (the latter
set must be suitably enlarged to match the size of the previous one
adding column indices greater than the indices
in $J \cup L$). Apply $S$ to the product of the minors:
$$
S([K,L][I,J])=[K_0,L_1][I_0,J_1] 
$$
where $K=K_0 \cup R$, $I=I_0 \cup R$ and $L_1 = \{1 \dots n \} \- J$,
$J_1 = \{1 \dots n \} \- L$.
Since $I_0 \cap K_0 = \emptyset$ we can apply the Theorem (5.5):
$$
\align
q^{-z}[K_0,L_1][I_0,J_1]
\equiv &
\sum_{i=0}^{M+N} (q^{-1}-q)^i \sum_{(Z,W) \in C_i^{I \- R,J \- C}}
sign(Z,W) (-q)^{-l(W)-l(W')} \cr\cr
&[Z,(W \cup C_0)_{ord}][Z',(W' \cup C_0)_{ord}]
\endalign
$$
where $C_0= J_1 \cap L_1$, $|C_0|=z$.
If we apply again the antiendomorphism $S$ we obtain the result.

\medskip
{\bf Example (7.4)}. Consider $[345,134]$, $[123,125]$. 
$R=\{3\}$, $C=\{1\}$. $I \cup K=\{12345\}$, $J \cup K=\{12345\}$. 
Apply $S$:
$$
S([345,134][123,125])=[45,34][12,25]
$$
We need to consider standard towers for row index $I_0=(12)$ in 
$I \cup K \- R= \{1245\}$ and for column index $(25)$ in
$J \cup L \- C= \{2345\}$. The only non trivial standard tower
is the column one:
$$
(25) <_{(23)} (35) <_{(34)} (45)
$$
The sets $C_i^{I \- R, J \- C}$ are:
$$
\matrix
s & \s_s& C_0^{(I \- R)_0, (J \- C)_s} & C_1^{(I \- R)_0, (J \- C)_s} & 
 C_2^{(I \- R)_0, (J \- C)_s}\cr\cr
0 & & (12,45)(45,23) &  &\cr\cr
1 & (34) & (12,35)(45,24) & (12,45)(45,23) & \cr\cr
2 & (23) & (12,25)(45,34) & (12,35)(45,24) & (12,45)(45,23) \NO \cr\cr
 &       &                & (12,45)(45,32)
\endmatrix
$$
Hence the commutation is
$$
\align
[345,134][123,125]= &
[123,125][345,134]+\cr\cr
&-(q^{-1}-q)[[123,135][345,124]-q^{-1}[123,145][345,123]]
\endalign
$$

\medskip
We want to make some remarks on the connection between the commutation
relations and two orderings.

\medskip

It is immediate to check that if $(K,L) > (I,J)$
lexicographically, then in (7.5) are listed all the propositions
that give in all possible cases the commutation relations
between the two minors $[K,L]$, $[I,J]$ according to the
definition (4.2). 

\medskip

It is also possible to refine this result. Consider the order appeared
in [GL]:
$$
(K,L) <_{\GL} (I,J) \qquad \hbox{iff} \qquad K \leq_r I, L \leq_c J 
$$
where $(k_1 \dots k_s) \leq_r (i_1 \dots i_r)$ if and only if
$s<r$ or $k_\a \leq i_\a$, for all $\a$'s. 
$(l_1 \dots l_s) \leq_c (j_1 \dots j_r)$ if and only if
$s< r$ or $j_\b \leq l_\b$, for all $\b$'s. 
\medskip

Then from the construction of a standard row and column tower
it is immediate that the commutation relation of the two 
minors $[K,L]$, $[I,J]$,
$(K,L) > (I,J)$ (lexicographically) 
is a product of minors $[\tilde I,\tilde J][\tilde K,\tilde L]$
such that $[\tilde I,\tilde J]<_{\GL}[I,J]$.

\medskip

{\bf Summary (7.5)}
Here is a summary of all the possible cases of commutation relations
among two quantum minors $[I,J]$, $[K,L]$.

\medskip

1. Case $I \cap K = \emptyset$, $J \cap L = \emptyset$

$intr(I)=0$, $intc(J)=0$: Lemma (4.4).

$intr(I)=0$, $intc(J) \geq 0$: Theorem (5.1).

$intr(I) \geq 0$, $intc(J)=0$: Theorem (7.1).

$intr(I)\geq 0$, $intc(J) \geq 0$: Theorem (6.1).

\medskip
2. Case $I \cap K = \emptyset$, $J \cap L \neq \emptyset$

$intr(I)=0$, $intc(J)=0$: Proposition (5.4).

$intr(I)=0$, $intc(J) \geq 0$: Theorem (5.5).

$intr(I) \geq 0$, $intc(J)=0$ and
$intr(I)\geq 0$, $intc(J) \geq 0$: Theorem (6.3).

\medskip

3. Case $I \cap K \neq \emptyset$, $J \cap L = \emptyset$

$intr(I)=0$, $intc(J)=0$ and $intr(I) \geq 0$, $intc(J)=0$: Theorem (7.1).

$intr(I)=0$, $intc(J) \geq 0$ and $intr(I)\geq 0$, $intc(J) \geq 0$: 
Theorem (7.2).

\medskip

4. Case $I \cap K \neq \emptyset$, $J \cap L \neq \emptyset$: Theorem (7.3).

\enddocument